\documentclass[12pt]{article}

\usepackage[T2A]{fontenc}
\usepackage[utf8]{inputenc}

\usepackage{amsfonts}
\usepackage{amssymb}
\usepackage{amsmath}
\usepackage{amsthm}

\def \le {\leqslant}
\def \ge {\geqslant}

\theoremstyle{plain}

\begin{document}
\begin{Large}
 \centerline{О критерии  Лежандра   подходящей дроби}
\end{Large}
\vskip+0.8cm
\begin{large}
\centerline{ Н.Г. Мощевитин, А.Ю. Яшникова }

\end{large}
\vskip+0.3cm
\centerline{ 
    (Астраханский государственный университет)}
\vskip+2.5cm

Как известно, вещественные числа представимы в виде обыкновенных цепных дробей
\begin{equation}\label{0}
   [a_0;a_1,a_2,...,a_n,...]
   =   a_0  +
\frac{1}{\displaystyle{a_1+\frac{1}{\displaystyle{a_2 + \cdots+
\frac{1}{\displaystyle{a_n + ...
{} }}}}}}
   , a_0\in \mathbb{Z},\,\,\,\ a_n\in \mathbb{Z}_+, n=1,2,3,... ,
\end{equation}
конечных для рациональных чисел и бесконечных для иррациональных, а  выражения
$$
\frac{p_n}{q_n} =    [a_0;a_1,a_2,...,a_n], \,\,\,\,\, (p_n,q_n) = 1,\,\,\,\,\,  n=1,2,3,... 
$$
называются подходящими дробями для дроби (\ref{0}).

В этой короткой заметке мы изложим любопытную на наш взгляд историю, связанную с известным  критерием того, что заданное рациональное число
$\frac{p}{q},  (p,q)= 1$  является подходящей дробью для некоторого  вещественного числа $\alpha$.  В особенности примечательным кажется нам вопрос, почему 
Люка  \cite{LU}  в процессе изложения этого критерия никоим образом не цитирует книгу  Лежандра  \cite{lege17}. Ответа на него мы не знаем.
В своей книге Люка упоминает имя Лежандра всего только около десяти раз, причем в главе XXIV, посвященной цепным дробям,  при обсуждении критерия подходящей дроби (параграф 245) Люка на  Лежанд-ра  нигде не ссылается.
Мы надеемся, что наш рассказ пробудит интерес к первоисточникам и поспособствует тому, чтобы читатель сам посмотрел и разобрался, что именно написано в книгах 
 \cite{lege17} и   \cite{LU}, 
 а также и в других упомянутых нами первоисточниках, которые сейчас общедоступны и легко могут быть найдены  в интернете.

  Обсуждаемая здесь история в какой-то мере  затронута в недавней статье первого автора \cite{Moshe}, однако сейчас все будет описано более подробно и обстоятельно.
  
  \vskip+0.3cm
  
 {\bf 1. Наиболее известная формулировка.}

  \vskip+0.3cm
  Во многих учебниках и книгах, где обсуждаются простейшие свойства цепных дробей,  
  следующее утверждение, которое мы здесь назовем теоремой А, формулируется как теорема Лежандра.

  В качестве примера, когда теорема А  формулируется как принадлежащая Лежандру, мы хотим привести следующие ссылки:  
  теорема 5С главы I из \cite{sch},
  теорема 1.8 из книги \cite{buj}.  Конечно, утверждение теоремы А очень хорошо известно, и  оно имеется в большинстве учебников 
  (см., например, \cite{X,bu,ne}).

---------------------------------------------------------------------------------

  \begin{small}
 Работа поддержана грантом РФФИ номер 18-01-00886.
  \end{small}

   \newpage
  
  {\bf Теорема А.}
   {\it  Расмотрим  вещественное число  $\alpha$. Если несократимая дробь $\frac{p}{q}$ удовлетворяет неравенству
 \begin{equation}\label{lll}
 \left|\alpha -\frac{p}{q}\right|<\frac{1}{2q^2}
 \end{equation}
то она является одной из подходящих дробей к числу $\alpha$.}
  
    \vskip+0.3cm

  Итак, теорема А дает достаточное условие того, что некоторая несократимая дробь $\frac{p}{q}$  являет-ся подходящей дробью для $\alpha$. Доказательство теоремы А зачастую излагается довольно причу-дливо (см., например \cite{sch}),  часто теорему А доказывают с помощью теоремы Лагран-жа, утверждающей, что  подходящие дроби и наилучшие приближения
  {\it второго рода } суть одно и то же (например, так делается в \cite{X}).
  На самом деле, теорема А является следствием более сильного  утверждения, доказанного Лежандром, а именно необходимого и достаточного условия  (то есть критерия) того, что дробь
  $\frac{p}{q}$  является подходящей дробью. Несколько  странно то, что, несмотря на то, что  точный критерий и его простое и  естественное доказательство 
  имеются уже в первом издании книги Лежандра  \cite{lege17} и, естественно, содержатся в последующих изданиях (например, \cite{lege}), критерий этот
  несколько менее известен, чем  утверждение теоремы А.

   \vskip+0.3cm
 {\bf  2.  Критерий Лежандра. } 
   \vskip+0.3cm
 
 Мы не будем здесь приводить дословную формулировку  критерия из   самой книги Лежандра  (см.  \cite{lege17}, p.27 - 29, или \cite{lege}, p.23 - 26, параграф (9)). В те времена не было принято всегда выделять теоремы отдельными утверждениями, и  Лежандр описывает содержание критерия
  в процессе его доказательства. 
  В параграфе (9) результат  теоремы А, которую мы сформулировали выше, он приводит в самом конце в последнем предложении.

   В книге Перрона  
  (мы  ссылаемся на 
   параграф 13 главы II из третьего издания  \cite{Per})
  критерий Лежандра сформулирован в виде замысловатого утверждения,
  из которого он  выводит Теорему А (Satz 2.11).  
  Мы решили не приводить дословно формулировку из книги Перрона тоже. Книга Лежандра написана по-французски,
  книга Перрона --- по-немецки,  и нам не хочется приводить здесь перевод утверждения. Но мы
   почти дословно  приведем формулировку из книги Венкова  (глава II,  параграф  5 из \cite{Ve},
  которая почти тождественна и  ориги-нальной формулировке  Лежандра, и   утверждению из книги Перрона, хотя в ней есть некоторый любопытный нюанс.
  
     \vskip+0.3cm
     
  {\bf Теорема Б (Критерий Лежандра).}
   {\it  Расмотрим  вещественное число  $\alpha$.  Считая  $\frac{p}{q}\neq \alpha$,
   положим
     \begin{equation}\label{1}
   \alpha - \frac{p}{q} = \frac{\theta}{q^2}
 \end{equation}
   и разложим  $\frac{p}{q}$ в непрерывную дробь с нечетным  или четным\footnote{
   Обычно нумерация неполных частных в цепной дроби именно такая, как в формуле (\ref{0}),
   в частности, так обстоит дело в книге Перрона. Но в книге Венкова нумерация неполных частных начинается с единицы, то есть там считается, что
   $[\alpha] = a_1$. Поэтому мы позволили себе изменить оригинальный текст формулировки теоремы Б и поменять местами слова "нечетным" и "четным".
   У Лежандра неполные частные  вообще не нумеруются индексами и обозначаются $\alpha, \beta, \gamma,\delta,...$.
   }
    числом неполных частных, смотря по тому, будет ли 
   $\theta<0$ или  
   $\theta>0$;
   пусть
   $\frac{p'}{q'}$ будет предпоследняя подходящая дробь к этой непрерывной дроби. Для того, чтобы $\frac{p}{q}$ была подходящей дробью к 
   $\alpha$,
   необходимо и достаточно, чтобы 
   \begin{equation}\label{m}
   |\theta| \le \frac{q}{q+q'}.
   \end{equation}
     }
     
          \vskip+0.3cm
          
     Различие между оригинальной формулировкой Лежандра (которая повторена у Перрона) и формулировкой теоремы Б  из книги Венкова состоит в том, что 
     у Лежандра и Перрона в качестве условия эквивалентного тому,  что $\frac{p}{q}$ 
     есть подходящая дробь к $\alpha$ вместо  {\it нестрогого} неравенства (\ref{m}) 
       фигурирует {\it строгое} неравенство
        \begin{equation}\label{1m}
          |\theta| < \frac{q}{q+q'}.
\end{equation}
  
  По-видимому, и сама формулировка теоремы Б, и  ее различие с тем, что имеется у Лежандра и Перрона, нуждаются в небольшом комментарии.
  
  Во-первых, всякое  рациональное число $\frac{p}{q}$ может быть  записано в виде обыкновенной цепной дроби ровно двумя способами
  \begin{equation}\label{raz}
  \frac{p}{q}=
     [a_0;a_1,a_2,...,a_{n-1},a_n] =
   [a_0;a_1,a_2,...,a_{n-1},a_n-1,1] 
  ,\,\,\,\,
   a_n \ge 2,
  \end{equation}
  и ясно, что здесь одно из разложений четной длины, а другое  - нечетной. Таким образом, мы пояснили, что имеется в виду в формулировке теоремы Б, и предпоследняя дробь  $\frac{p'}{q'}$  будет, соответственно, определяться равенствами 
   \begin{equation}\label{raa}
   \frac{p'}{q'} =
     [a_0;a_1,a_2,...,a_{n-1}] \,\,\,\,\,
    {или}
       \frac{p'}{q'} =
   [a_0;a_1,a_2,...,a_{n-1},a_n-1] 
 .
  \end{equation}
  
  Во-вторых,  если для $\theta$, определенного в 
  (\ref{1}),  выполнено равенство
   $$
          |\theta|=\frac{q}{q+q'},
  $$
  то это значит, что 
   $$
   \alpha
   = \frac{p+p'}{q+q'}
   $$
   является медиантой дробей
  $$
  \frac{p}{q}=
     [b_0;b_1,b_2,...,b_{k-1},b_k]
     \,\,\,\,\,
     \text{и}
     \,\,\,\,\,
    \frac{p'}{q'}=
     [b_0;b_1,b_2,...,b_{k-1}]      
  ,\,\,\,\,\,
   k = n \,\,\, \text{или}\,\,\, n+1.
  $$
  В этом случае $ \alpha$ рационально, и для него имеется два разложения в цепную дробь:
  \begin{equation}\label{a}
  \alpha =    [b_0;b_1,b_2,...,b_{k-1},b_k,1]
  \end{equation}
  и
   \begin{equation}\label{a1}
  \alpha =    [b_0;b_1,b_2,...,b_{k-1},b_k+1].
  \end{equation}
  Для разложения (\ref{a})  дробь $\frac{p}{q}$ будет (предпоследней) подходящей дробью, но для разложения
   (\ref{a1})  дробь $\frac{p}{q}$   подходящей дробью не будет. Таким образом, свойство дроби $\frac{p}{q}$  быть подходящей дробью для  $ \alpha$ зависит от того, какое разложение мы взяли, и  чтобы пояснить различие,  снова надо апеллировать к неединственности разложения рационального числа в цепную дробь.

          \vskip+0.3cm

  {\bf 3. Формулировка Люка.}

          \vskip+0.3cm

          Утверждение, похожее на критерий Лежандра, приведено в книге Люка \cite{LU}. Мы уже отмечали в начале нашей заметки, что Люка почему-то не ссылается на Лежандра. Более того, формулировка Люка не совсем четкая. Постаравшись быть близкими к французскому оригиналу,  мы приводим ее в виде нижеследующей теоремы В.
          
             \vskip+0.3cm
  
  {\bf Теорема В.}
   {\it   Для того, чтобы данная дробь $f\!:\!g$ 
   была бы подходящей дробью к числу $x$, необходимо и достаточно выполнение неравенства
   \begin{equation}\label{x1}
   \left| x -\frac{f_p}{g_p}\right|<
   \frac{1}{g_p(g_p+g_{p-1})},
   \end{equation}
   где через 
   $g_p $ и $ g_{p-1}$
   обозначены знаменатели двух последних подходящих дробей разложения  величины $f\!:\!g$ 
   в цепную дробь.

   }
  
   \vskip+0.3cm
   
   Заметим, что Люка упоминает, что в предположении (\ref{lll}) неравенство (\ref{x1})  для $x=\alpha$ выполня-ется.

  Теперь мы обратим наше внимание на следующее утверждение, которое использовалось в статье Лемера \cite{lehmer}. Мы приводим дословный перевод формулировки  оттуда:
  
   \vskip+0.3cm
  
  {\bf Теорема Г.}
   {\it 
   Пусть $A_k/B_k, A_{k-1}/B_{k-1}$ суть две последовательные  подходящие дроби к числу $\eta$. Тогда эти две дроби будут двумя последовательными подходящими дробями для числа $\xi$ тогда и только тогда, когда
   \begin{equation}\label{x}
   \left| \xi - \frac{A_k}{B_k}\right|
   < \frac{1}{B_k(B_k+B_{k-1})}.   
\end{equation}
   }
  
   \vskip+0.3cm

  Это утверждение Лемер называет хорошо известным и дает ссылку как раз на книгу Люка \cite{LU}.   
  Формулировка Лемера несколько отличается от формулировки Люка, но она тоже не совсем четкая.
  Ниже мы прокомментируем  утверждения теорем В и Г.
  
  Начнем с утверждение Люка - теоремы В. На первый взгляд, неравенство (\ref{x1}) задает интервал с серединой $\frac{f_p}{g_p}$. Но на самом деле, это неравенство надо интерпретировать по-другому. Как мы уже писали выше, 
  рациональное число раскладывается в цепную дробь ровно двумя способами (\ref{raz}). В утверждении теоремы В надо
  брать два последних знаменателя разложения величины  $f\!:\!g$  в цепную дробь. Но имеется два различных разложения.
  И в зависимости от знака  разности $ x -\frac{f}{g}$  набо выбирать разложение четной или нечетной длины. Точнее, если
  $ x >\frac{f}{g}$, надо брать разложение  $  \frac{f}{g} =\frac{f_p}{g_p}$   с четным $p$, а если 
    $ x <\frac{f}{g}$,  то надо брать разложение    с нечетным $p$.

  Теперь прокомментируем утверждение теоремы Г.
  
    Во-первых, 
    число $\xi$ обязано лежать именно между соседними подходящими дробями  $A_k/B_k$  и $A_{k-1}/B_{k-1}$, поэтому теорему Г следует уточнить:  необходимым и достаточным условием является принадлежность  $\xi$  интервалу (\ref{x}), пересеченному с отрезком
    с концами  $A_k/B_k, A_{k-1}/B_{k-1}$.

Во-вторых, две несократимые дроби $\frac{p}{q}, \frac{p'}{q'}$ будут последовательными подходящими дробями к некоторому 
  $\eta$ тогда и только тогда, когда
  \begin{equation}\label{fraa}
  \left|
  \frac{p}{q}- \frac{p'}{q'}
  \right| = \frac{1}{qq'}.
  \end{equation}
  Следовательно, мы можем переформулировать теорему Г так:

     \vskip+0.3cm
  
  {\bf Теорема Г$'$.}
  {\it
  Две дроби $\frac{p}{q}, \frac{p'}{q'}$ с условием $ q\ge q'$ и удовлетворяющие
  (\ref{fraa}) 
  будут двумя последовательными подходящими дробями для числа $\xi$ тогда и только тогда, когда $\xi$ 
  принадле-жит полуоткрытому интервалу с концами    $\frac{p}{q}$ и $\frac{p+p'}{q+q'} =  \frac{p}{q}\pm \frac{1}{qq'}$,
  причем конец  $\frac{p}{q}$  входит в полуинтервал,  а  конец   $\frac{p+p'}{q+q'}  $ не входит.
  }
    \vskip+0.3cm

  \newpage
{\bf  4.  Более естественная формулировка.} 

 \vskip+0.3cm
 
Очень уместным нам кажется  сформулировать теорему Лежандра с помощью рядов Фарея. 
 Мы ограничимся случаем  $\frac{p}{q}\in (0,1)$.
Для заданной несократимой рациональной дроби  $\frac{p}{q}$ рассмотрим ряд Фарея $\frak{F}_q$ порядка  $q$, то есть совокупность всех
рациональных дробей со знамена-телем $\le q$:
$$
\frak{F}_q: \,\,\,\,\,0 = r_0< r_1<r_2<...<r_j< r_{j+1}<...<r_\Phi,\,\,\,\
\Phi = \sum_{k=1}^q \varphi(q),\,\,\,\,
\varphi(\cdot) - \text{функция Эйлера}.
$$
Где-то в этой последовательности находится число $\frac{p}{q} = r_j$.
 Мы рассмотрим  его двух соседей по ряду Фарея $\frak{F}_q$: 
 $$r_{j-1}< \frac{p}{q} = r_j< r_{j+1}.
 $$ 
 Если рассмотреть единственное разложение дроби  $\frac{p}{q}$ в цепную дробь вида
$$
r_j=
\frac{p}{q} = [a_0;a_1,..., a_t]
,\,\, a_t\ge 2,
$$
то для соседних дробей будет выполнено
$$
r_{j-1}= \frac{p_-}{q_-} = [a_0;a_1,...,a_{t-1}] ,\,\,\,
r_{j+1} =  \frac{p_+}{q_+}=[a_0;a_1,...,a_{t}-1] ,
$$
{или}
$$
r_{j-1}= \frac{p_-}{q_-} = [a_0;a_1,...,a_{t}-1] ,\,\,\,
r_{j+1} =  \frac{p_+}{q_+}=[a_0;a_1,...,a_{t-1}] ,
$$
в зависимости от четности $t$. Ясно, что теоремы Б  
(со строгим неравенством (\ref{1m})  вместо  (\ref{m})) и теорема В 
(уточненная нашим комментарием о том, как выбирать разложение в цепную дробь)
эквивалентны  следующему утверждению.

     \vskip+0.3cm
  {\bf Теорема Д.}
   {\it  
Несократимая дробь  $\frac{p}{q} $ будет подходящей дробью для числа
 $\alpha$ тогда и только тогда, когда
 $$
 \alpha \in \left(
 \frac{p+p_-}{q+q_-}, 
 \frac{p+p_+}{q+q_+}
 \right).
  $$}
     \vskip+0.3cm

 %Примечательно, что в книге Люка имеется параграф про последовательности Штерна-Броко.

 Кроме того, нам кажется естественным упомянуть здесь рассуждения (в духе параграфа  12 главы III книги Хинчина \cite{X}),    связанные с "цилиндрами" 
 \begin{equation}\label{cc}
 E
 \left(
 \begin{array}{ccccc}
 1 & 2& ... & n-1 & n
 \cr
  a_1 & a_2&  ... & a_{n-1} & a_n
 \end{array}
 \right) =
 \end{equation}
 $$
 =
 \{ \alpha \in [0,1]:
 \,
 \alpha \,
 \text{представляется в виде конечной или бесконечной дроби}\, 
 \alpha = [0; a_1,a_2,...,a_{n-1}, a_n,...]\}.
 $$
 Легко видеть что цилиндр  (\ref{cc})  представляет из себя полуинтервал  с концами
 $ \frac{p_n}{q_n}, \frac{p_n+p_{n-1}}{q_n + q_{n-1}}$,
 точнее
   $$
 E
 \left(
 \begin{array}{ccccc}
 1 & 2& ... & n-1 & n
 \cr
  a_1 & a_2&  ... & a_{n-1} & a_n
 \end{array}
 \right) =
 \begin{cases}
 \left[
  \frac{p_n}{q_n}, \frac{p_n+p_{n-1}}{q_n + q_{n-1}}\right),\,\, \text{если  $n$ четно},
\cr
  \left(
   \frac{p_n+p_{n-1}}{q_n + q_{n-1}},  \frac{p_n}{q_n}\right],\,\,
  \text{если   $n$ нечетно}
 .
 \end{cases}
 $$
 Именно этот полуинтервал   фигурирует в формулировке теоремы Г$'$.
 Как мы отмечали ранее,
 рациональное число $ \frac{p}{q}$ представляется в 
 виде цепной дроби ровно двумя способами  (\ref{raz}). Таким образом,
 $$
 \left(
 \frac{p+p_-}{q+q_-}, 
 \frac{p+p_+}{q+q_+}\right)=
  E
 \left(
 \begin{array}{ccccc}
 1 & 2& ... & n-1 & n
 \cr
  a_1 & a_2&  ... & a_{n-1} & a_n
 \end{array}
 \right) 
 \cup
  E
 \left(
 \begin{array}{cccccc}
 1 & 2& ... & n-1 & n&n+1
 \cr
  a_1 & a_2&  ... & a_{n-1} & a_n-1& 1
 \end{array}
 \right) ,
 $$
 причем, естественно,
 $$
   E
 \left(
 \begin{array}{ccccc}
 1 & 2& ... & n-1 & n
 \cr
  a_1 & a_2&  ... & a_{n-1} & a_n
 \end{array}
 \right) 
 \cap
  E
 \left(
 \begin{array}{cccccc}
 1 & 2& ... & n-1 & n&n+1
 \cr
  a_1 & a_2&  ... & a_{n-1} & a_n-1& 1
 \end{array}
 \right) =\left\{ \frac{p}{q}\right\}.
 $$
 
 В силу определения множеств $E(\cdot)$,
 последние равенства доказывают теорему Д.
% Впрочем, по сути дела, оригинальное доказательство Лежандра  устроено подобным же образом.

     \vskip+0.3cm

  {\bf 5. Об обобщениях.} 
  
      \vskip+0.3cm
  
  Следующее обобщение теоремы А со ссылками на Фату \cite{f}  и Грейса \cite{g}  приводится в книге  Бюжо \cite{buj}. В этой книге следующий результат назван "малоиз-вестным" (\cite{buj}, cтр. 9), см., однако, теорему 10 главы I из книги Ленга \cite{LA}.

      \vskip+0.3cm

  {\bf Теорема Е.} {\it 
 Если несократимая дробь $\frac{p}{q}$ удовлетворяет неравенству
 $$\left|\alpha -\frac{p}{q}\right|<\frac{1}{q^2}$$ 
то она является одной из  
подходящих или промежуточных дробей для $\alpha$, из набора
$$
\frac{p_n}{q_n}, \,\,\,\,\,
\frac{p_{n+1}+ p_n}{q_{n+1}+ q_n},\,\,\,\,\,
\frac{p_{n+1}- p_{n}}{q_{n+1}- q_{n}}
$$
при некотором 
$n$.
}
     \vskip+0.3cm

 Примечательно, что относительно недавно были опубликованы утверждения, содержащие еще более широкие простые  обобщения.
 Например, Дуелла \cite{duj} пишет, что он обобща-ет классический результат Лежандра - теорему А и доказывает следующее утверждение.

     \vskip+0.3cm
   {\bf Теорема Ж.} {\it 
   Пусть $\alpha$ - иррациональное число и несократимая дробь $\frac{p}{q}$ удовлетворяет неравенству
   $$\left|\alpha -\frac{p}{q}\right|<\frac{c}{q^2}$$ 
   с некоторым положительным 
   $c$.
   Тогда пара $(p,q)$  имеет вид
   $$
   (p,q) = (rp_{n+1}\pm sp_n, rq_{n+1}\pm sq_n)
   $$
   с некоторыми неотрицательными целыми $r,s$,
   удовлетворяющими   неравенству $rs <2c$.}
   
       \vskip+0.3cm
   
   Несколько более громоздкая формулировка имеется  более ранней работе Ворли
   \cite{wo}.

  \end{document}